\documentstyle[12pt]{article}
\begin{document}\thispagestyle{empty}\begin{flushright}
OUT--4102--71\\
math.CA/9803067\\
16 March 1998          \end{flushright}\vspace*{2mm}\begin{center}{\Large\bf
Polylogarithmic ladders,
hypergeometric series                              \\[5pt]
and the ten millionth digits
of $\zeta(3)$ and $\zeta(5)$                       }\vglue 10mm{\large{\bf
D.~J.~Broadhurst                                   $^{1)}$}\vglue 4mm
Physics Department, Open University                \\[3pt]
Milton Keynes MK7 6AA, UK     }\end{center}\vfill\noindent{\bf Abstract}\quad
We develop ladders that reduce $\zeta(n):=\sum_{k>0}k^{-n}$, for
$n=3,5,7,9,11$, and $\beta(n):=\sum_{k\ge0}(-1)^k(2k+1)^{-n}$,
for $n=2,4,6$, to convergent polylogarithms and products of powers of $\pi$
and $\log2$. Rapid computability results because the required arguments of
${\rm Li}_n(z)=\sum_{k>0}z^k/k^n$ satisfy $z^8=1/16^p$, with  $p=1,3,5$.
We prove that $G:=\beta(2)$, $\pi^3$, $\log^32$, $\zeta(3)$, $\pi^4$,
$\log^42$, $\log^52$, $\zeta(5)$, and six products of powers of $\pi$ and
$\log2$ are  constants whose $d$th hexadecimal digit can be computed in
time~$=O(d\log^3d)$ and space~$=O(\log d)$, as was shown for $\pi$, $\log2$,
$\pi^2$ and $\log^22$ by Bailey, Borwein and Plouffe. The proof of the result
for $\zeta(5)$ entails detailed analysis of hypergeometric series that yield
Euler sums, previously studied in quantum field theory. The other 13 results
follow more easily from Kummer's functional identities. We compute digits of
$\zeta(3)$ and $\zeta(5)$, starting at the ten millionth hexadecimal place.
These constants result from calculations of massless Feynman diagrams
in quantum chromodynamics. In a related paper, hep-th/9803091, we
show that massive diagrams also entail constants whose base of super-fast
computation is $b=3$.
\vfill\footnoterule\noindent
$^1$) D.Broadhurst@open.ac.uk;
http://physics.open.ac.uk/$\;\widetilde{}$dbroadhu
\newpage\setcounter{page}{1}
\newcommand{\dfrac}[2]{\mbox{$\frac{#1}{#2}$}}

\section{Introduction}

David Bailey, Peter Borwein and Simon Plouffe~\cite{BBP} showed that
$\pi$, $\log2$, $\pi^2$ and $\log^22$
are in the class SC$^*$ of constants whose $d$th digit can be computed in
time~$=O(d\log^{O(1)}d)$ and space~$=O(\log^{O(1)}d)$.
They gave simple methods to find the $d$th hexadecimal
digit of each of these constants in time~$=O(d\log^3d)$ and
space~$=O(\log d)$. Subsequently, Fabrice Bellard
found a slightly faster~\cite{FB} algorithm for $\pi$
and computed it at $d=250,000,000,000$. It was not known if Catalan's
constant, $\pi^3$, $\log^32$, $\zeta(3)$, $\pi^4$, $\log^42$, $\log^52$,
or $\zeta(5)$, are in this class.
We shall prove that they are, by studying polylogarithms~\cite{Lew}
${\rm Li}_{n}(z)=\sum_{k>0} z^k/k^n$ with $z^8=1/16^p$ and
$p=1,3,5$.

Section~2 proves all but one of the required identities,
staying within the realm of functional relations of polylogarithms.
To complete the proof for $\zeta(5)$, we derive a hypergeometric
generating function in section~3, where the process is also reversed,
to deduce further hypergeometric identities from the results that
were proved by functional relations. Digits of $\zeta(3)$ and $\zeta(5)$
are given in section~4.

\section{Polylogarithmic ladders to order $n=11$}

The 8 polylogarithmic ladders that figure here originate in rather simple
observations
\begin{eqnarray}
w:=\frac{1+i}{2}&\Longrightarrow&w=\frac{1+w^3}{1-w^2}
=\frac{1+w^5}{(1-w^2)^2}\label{w}\\
h:=\frac{i}{\sqrt2}&\Longrightarrow&h=\frac{1-h^3}{(1-h)^2}\label{h}\,.
\end{eqnarray}

\subsection{Logarithms}

In terms of ${\rm Li}_1(z):=-\log(1-z)=\sum_{k>0}z^k/k$,
observations~(\ref{w},\ref{h}) give
\begin{eqnarray}
{\rm Li}_1(w)-
\dfrac12\,{\rm Li}_1(\dfrac12)&=&\phantom{+}i\pi/4\label{w11}\\
{\rm Li}_1(-w^3)-{\rm Li}_1(w^2)-
\dfrac12\,{\rm Li}_1(\dfrac12)&=&-i\pi/4\label{w13}\\
{\rm Li}_1(-w^5)-2\,{\rm Li}_1(w^2)-
\dfrac12\,{\rm Li}_1(\dfrac12)&=&-i\pi/4\label{w15}\\
{\rm Li}_1(h^3)-2\,{\rm Li}_1(h)-
\dfrac12\,{\rm Li}_1(\dfrac12)&=&-i\pi/2\label{h1}
\end{eqnarray}
which provide a variety of ways of evaluating $\pi$.
For example, one may multiply~(\ref{w11}) by $8(1-w)=4(1-i)$
and take the real part, to obtain
\begin{eqnarray}
\pi&=&
8\,\Re\left\{(1-w)\,{\rm Li}_1(w)\right\}-2\,{\rm Li}_1(\dfrac12)\nonumber\\
&=&\sum_{k\ge0}\frac{1}{16^k}\left\{\frac{4}{8k+1}-\frac{2}{8k+4}
-\frac{1}{8k+5}-\frac{1}{8k+6}\right\}\label{BBP}
\end{eqnarray}
without need of the integration used to prove this
result in~\cite{BBP}.

The features of note here are
relations~(\ref{w13}--\ref{h1}), with powers $p=2,3,5$.
We shall prove polylogarithmic generalizations that yield 14 new SC$^*$
results. The ladders that produce these results
will be extended, in sections~2.2--2.11, to polylogarithms of order $n=11$.

\subsection{Dilogarithms}

Rewriting the order $n=1$ relations~(\ref{w11}--\ref{w15}) as
\begin{eqnarray}
{\rm Li}_1(w)&=&-\log(1-w)\label{w11b}\\
{\rm Li}_1(-w^3)&=&{\rm Li}_1(w^2)-\log(w)\label{w13b}\\
{\rm Li}_1(-w^5)&=&2\,{\rm Li}_1(w^2)-\log(w)\label{w15b}
\end{eqnarray}
we find $n=2$ relations that are strikingly similar, namely
\begin{eqnarray}
2\,{\rm Li}_2(w)&=&-\log^2(1-w)-2\,{\rm Li}_2(-i)\label{w21}\\
2\,{\rm Li}_2(-w^3)&=&3\left\{{\rm Li}_2(w^2)-\log^2(w)\right\}
+4\,{\rm Li}_2(-i)\label{w23}\\
2\,{\rm Li}_2(-w^5)&=&5\left\{2\,{\rm Li}_2(w^2)-\log^2(w)\right\}
+8\,{\rm Li}_2(-i)\label{w25}
\end{eqnarray}
with ${\rm Li}_2(-i)=-i G-\pi^2/48$, where $G$ is Catalan's constant.

The method used to prove~(\ref{w21}--\ref{w25})
was designed to capture all identities with powers $p\le5$. It was
also designed to be generalizable to orders $n=3,4,5$.
We took Kummer's functional relation for {\rm Li}$_2(x(1-y)^2/y(1-x)^2)$,
given in Eq~(A.2.1.19) of Leonard Lewin's
book~\cite{Lew}, and substituted $x$ and $y$ by values drawn from
$\{\frac12,-1,i,w,1-w\}$ and their inverses.
After using duplication and inversion
relations that apply at any order, one sees that all
such specializations of Kummer's functional relations generate
polylogarithms with arguments restricted to $\{\frac12,-\frac12,-\frac18,
i,w,w^2,-w^3,-w^5\}$ and their complex conjugates. Thus one merely
has to find all the relations between 13 terms implied by $9^2=81$ equations,
which was a task easily performed by {\sc reduce}~\cite{REDUCE}
at each order. At $n=2$ this automated method proves~(\ref{w21}--\ref{w25})
and the real results
\begin{eqnarray}
\Re\,{\rm Li}_2(i)&=&-\dfrac{1}{48}\pi^2\label{h21}\\
{\rm Li}_2(\dfrac12)&=&\dfrac{1}{12}\pi^2-\dfrac12\log^22\label{h22}\\
\Re\left\{{\rm Li}_2(h^3)-6\,{\rm Li}_2(h)\right\}&=&
\dfrac{1}{12}\pi^2-\dfrac{3}{8}\log^22\,.\label{h23}
\end{eqnarray}

We remark that in~\cite{BBP} it was not possible to reduce
Catalan's constant to SC$^*$ sums with power $p=1$,
because~(\ref{w21}) determines only a combination of $G$ and $\pi\log2$.
The remedy is now to hand, in the $p\le3$ result~(\ref{w23}),
which furnishes another combination.
For brevity's sake, we define the rapidly computable series
\begin{equation}
S_{n,p}(a_1\ldots a_8):=\sum_{k>0}\frac{a_k}{2^{\lfloor
\frac{pk+p}{2}\rfloor}k^n}\label{S}
\end{equation}
on the understanding that the integer sequence $a_k$ has period 8 and is
hence specified by the arguments $a_1\ldots a_8$. Then the
results of~\cite{BBP} and~\cite{FB}, obtained from~(\ref{w11})
and~(\ref{w15}), are
\begin{eqnarray}
\pi&=&8\,S_{1,1}(1,0,0,-1,-1,-1,0,0)\label{pi}\\
&=&16\,S_{1,1}(0,1,0,0,0,-1,0,0)-16\,S_{1,5}(1,1,1,0,-1,-1,-1,0)\,.\label{FB}
\end{eqnarray}
At $n=2$ and $p=1$, one obtains from~(\ref{w21},\ref{h22}) the results
of~\cite{BP}
\begin{eqnarray}
\pi^2&=&
32\,S_{2,1}(1,-1,-1,-2,-1,-1,1,0)\label{pi2}\\
\log^22&=&\dfrac83\,S_{2,1}(2,-5,-2,-7,-2,-5,2,-3)\,.\label{l2}
\end{eqnarray}
Our new result for Catalan's constant, proved by~(\ref{w21},\ref{w23}), is
\begin{equation}
G=3\,S_{2,1}(1,-1,1,0,-1,1,-1,0)
-2\,S_{2,3}(1,1,1,0,-1,-1,-1,0)\,.\label{G}
\end{equation}

Some obvious questions arise. When do
\begin{eqnarray}
\lambda(n)&:=&\sum_{k\ge0}(2k+1)^{-n}=(1-2^{-n})\zeta(n)\label{lamn}\\
\beta(n)&:=&\sum_{k\ge0}(-1)^k(2k+1)^{-n}\label{betn}
\end{eqnarray}
cease to be in SC$^*$? We know that $\lambda(2k)$ and $\beta(2k+1)$
are rational multiples of $\pi^{2k}$ and $\pi^{2k+1}$. When
do {\bf Q}--linear combinations of $\lambda(2k+1)$ and $\beta(2k)$
with products of powers of $\pi$ and $\log2$
cease to be in SC$^*$?

To pursue these issues, one should pay heed to the integers
in~(\ref{w23},\ref{w25}).
We know from~\cite{BBP,Lew2}
how to extend to order $n=5$ a ladder containing~(\ref{h23}),
by inclusion of increasing powers of 2 and 3
at increasing orders. Supposing that~(\ref{w23})
shows the beginning of a similar pattern, and noting the first power of 5
in~(\ref{w25}), we define 8 ladders by
\begin{eqnarray}
A_n&:=&{\rm Li}_n(\dfrac12)\label{an}\\
B_n&:=&\Re\left\{2^{n-1}{\rm Li}_n(\dfrac{1+i}{2})\right\}\label{bn}\\
C_n&:=&\Re\left\{(2/3)^{n-1}{\rm Li}_n(\dfrac{i}{\sqrt8})
             -2^n{\rm Li}_n(\dfrac{-i}{\sqrt2})\right\}\label{cn}\\
D_n&:=&\Re\left\{(2/3)^{n-1}{\rm Li}_n(\dfrac{1+i}{4})
             -{\rm Li}_n(\dfrac{-i}{2})\right\}\label{dn}\\
E_n&:=&\Re\left\{(2/5)^{n-1}{\rm Li}_n(\dfrac{1-i}{8})
               -2\,{\rm Li}_n(\dfrac{-i}{2})\right\}\label{en}\\
F_n&:=&\Im\left\{2^{n-1}{\rm Li}_n(\dfrac{1+i}{2})\right\}\label{fn}\\
G_n&:=&\Im\left\{(2/3)^{n-1}{\rm Li}_n(\dfrac{1+i}{4})
             -{\rm Li}_n(\dfrac{-i}{2})\right\}\label{gn}\\
H_n&:=&\Im\left\{(2/5)^{n-1}{\rm Li}_n(\dfrac{1-i}{8})
               -2\,{\rm Li}_n(\dfrac{-i}{2})\right\}\label{hn}
\end{eqnarray}
and powers of $\log2$ by
\begin{equation}
L_n:=(-\log2)^n/n!\label{ln}
\end{equation}
with the implication that $L_n=0$ when $n<0$.
We suppress the appearance of $\log2$ and $\pi$, at $n=1,2$, by forming
the combinations
\begin{eqnarray}
\overline{A}_n&:=& A_n        +L_n-\dfrac{1}{2} \zeta(2)L_{n-2}\label{ab}\\
\overline{B}_n&:=& B_n+\dfrac12L_n-\dfrac{5}{8} \zeta(2)L_{n-2}\label{bb}\\
\overline{C}_n&:=& C_n+\dfrac12L_n-\dfrac{1}{3} \zeta(2)L_{n-2}\label{cb}\\
\overline{D}_n&:=& D_n+\dfrac12L_n-\dfrac{5}{24}\zeta(2)L_{n-2}\label{db}\\
\overline{E}_n&:=& E_n+\dfrac12L_n-\dfrac{7}{40}\zeta(2)L_{n-2}\label{eb}\\
\overline{F}_n&:=& F_n-\beta(1)L_{n-1}\label{fb}\\
\overline{G}_n&:=& G_n-\beta(1)L_{n-1}\label{gb}\\
\overline{H}_n&:=& H_n-\beta(1)L_{n-1}\label{hb}
\end{eqnarray}
with $\beta(1)=\pi/4$ and $\zeta(2)=\pi^2/6$. Then the results for $n=1,2$
are
\begin{eqnarray}
\overline{A}_1=\overline{B}_1=\overline{C}_1=\overline{D}_1=\overline{E}_1
=\overline{F}_1=\overline{G}_1=\overline{H}_1&=&0\label{r1}\\
\overline{A}_2=\overline{B}_2=\overline{C}_2=\overline{D}_2=\overline{E}_2
&=&0\label{r2}\\
\dfrac12\overline{F}_2=\dfrac34\overline{G}_2=\dfrac58\overline{H}_2&=&G:=
\beta(2)\,.\label{i2}
\end{eqnarray}

\subsection{Trilogarithms}

The machinery of proof employed at $n=2$ is now used at $n=3$,
by making substitutions in the appropriate Kummer
identity, which is Eq~(A.2.6.11) of~\cite{Lew}.
We obtain 8 relations between the 13 unknowns, namely
\begin{eqnarray}
\lambda(3)&=&\overline{A}_3=\dfrac25\overline{B}_3=\dfrac97\overline{C}_3
=3\overline{D}_3=\dfrac{25}{6}\overline{E}_3\label{r3}\\
\beta(3)&=&\dfrac{1}{32}\pi^3=\dfrac23\overline{F}_3-\overline{G}_3
=\dfrac{20}{23}\overline{F}_3-\dfrac{25}{23}\overline{H}_3\label{i3}
\end{eqnarray}
from which it follows that $\{\zeta(3),\pi^3,\pi^2\log2,\pi\log^22,
\log^32\}$ are in SC$^*$. In particular
\begin{eqnarray}
\log^32&=&192\,S_{3,1}(0,1,0,4,0,1,0,16)\nonumber\\
&&{}-32\,S_{3,3}(4,-3,-4,-1,-4,-3,4,7)\label{l3}\\
\lambda(3)=\dfrac78\zeta(3)&=&
6\,S_{3,1}(1,-7,-1,10,-1,-7,1,0)\nonumber\\
&&{}+4\,S_{3,3}(1,1,-1,-2,-1,1,1,0)\label{z3}\\
\beta(3)=\dfrac{1}{32}\pi^3&=&
5\,S_{3,1}(1,-6,1,0,-1,6,-1,0)\nonumber\\
&&{}+\dfrac53\,S_{3,3}(1,1,1,0,-1,-1,-1,0)\nonumber\\
&&{}+2\,S_{3,5}(1,1,1,0,-1,-1,-1,0)\label{b3}
\end{eqnarray}
whereas~\cite{BBP} found only two independent SC$^*$ combinations of
$\{\zeta(3),\log^32,\pi^2\log2\}$.

\subsection{Polylogarithms of order 4}

The expectations at $n=4$ are rather clear.
We expect to obtain $\pi^4$ from 4 combinations of $\overline{A}_4\ldots
\overline{E}_4$, and $\beta(4)$ from 2 combinations of $\overline{F}_4,
\overline{G}_4, \overline{H}_4$, and $\beta(3)\log2$. The method is as
before, taking all instances of Kummer's functional relation
that produce the target set of polylogarithms. At $n=4$ the relation is
Eq~(A.2.7.41) of~\cite{Lew}, which expresses $\log^2(1-x)\log^2(1-y)$
as a combination of 20 polylogarithms. The automated proof method
produces 7 relations. One gives $\lambda(4)=\pi^4/96$.
The other 6 are
\begin{eqnarray}
\overline{B}_4-\dfrac52\overline{A}_4&=&
\dfrac{343}{128}\zeta(4)\label{r4b}\\
\overline{C}_4-\dfrac79\overline{A}_4&=&
\dfrac{5}{54}\zeta(4)\label{r4c}\\
\overline{D}_4-\dfrac13\overline{A}_4&=&
-\dfrac{313}{3456}\zeta(4)\label{r4d}\\
\overline{E}_4-\dfrac{6}{25}\overline{A}_4&=&
-\dfrac{1547}{16000}\zeta(4)\label{r4e}\\
\overline{G}_4-\dfrac{2}{3}\overline{F}_4-\beta(3)\log2
&=&-\dfrac{80}{27}\beta(4)\label{i4g}\\
\overline{H}_4-\dfrac{4}{5}\overline{F}_4-\dfrac{23}{25}\beta(3)\log2
&=&-\dfrac{384}{125}\beta(4)\label{i4h}
\end{eqnarray}
with combinations on the left that could have been surmised
from the $n=3$ results. Note that two independent combinations
of $\beta(4)$, $\pi^3\log2$, and $\pi\log^32$ are in SC$^*$, but there
is no indication that $\beta(4)$ is. On the other hand, $\pi^4$,
$\pi^2\log^22$, and $\log^42$ are all obtained in SC$^*$, by choosing
3 of the 4 relations~(\ref{r4b}--\ref{r4e}). In particular, the first 3 give
\begin{eqnarray}
\dfrac{615}{256}\log^42&=&
3\,S_{4,1}(73,-2617,-73,-5066,-73,-2617,73,-27564)\nonumber\\
&&{}+S_{4,3}(1258,-761,-1258,-497,-1258,-761,1258,2019)
\label{l4}\\
\dfrac{41}{9216}\pi^4&=&
3\,S_{4,1}(1,-19,-1,-2,-1,-19,1,-108)\nonumber\\
&&{}+2\,S_{4,3}(3,-1,-3,-2,-3,-1,3,4)\,.\label{z4}
\end{eqnarray}

\subsection{Polylogarithms of order 5}

It is now proven that each of the 5 combinations
\begin{eqnarray}
\widetilde{B}_n&:=&
\overline{B}_n-\dfrac52\overline{A}_n
-\dfrac{343}{128}\zeta(4)L_{n-4}\label{bt}\\
\widetilde{C}_n&:=&\overline{C}_n-\dfrac79\overline{A}_n
-\dfrac{5}{54}\zeta(4)L_{n-4}\label{ct}\\
\widetilde{D}_n&:=&\overline{D}_n-\dfrac13\overline{A}_n
+\dfrac{313}{3456}\zeta(4)L_{n-4}\label{dt}\\
\widetilde{E}_n&:=&\overline{E}_n-\dfrac{6}{25}\overline{A}_n
+\dfrac{1547}{16000}\zeta(4)L_{n-4}\label{et}\\
\widetilde{H}_n&:=&
\overline{H}_n-\dfrac{4}{5}\overline{F}_n+\dfrac{23}{25}\beta(3)L_{n-3}
-\dfrac{648}{625}\left\{
\overline{G}_n-\dfrac{2}{3}\overline{F}_n+\beta(3)L_{n-3}\right\}\label{ht}
\end{eqnarray}
vanishes at $n=1,2,3,4$.
At $n=5$, we expect each of the first 4
to be rational multiples of $\lambda(5)$
and the last to be a rational multiple of $\beta(5)$.

The requisite functional identity for Li$_5$ is not given explicitly
in~\cite{Lew}. Rather, a combination of 34 polylogarithms is specified
that evaluates to products of logarithms and $\pi^2$.
The latter may be obtained by using 9 instances of
Kummer's 21-term functional relation at $n=4$.
Using {\sc reduce}~\cite{REDUCE} to perform this task,
we obtained the right--hand side of
\begin{eqnarray}&&
 {\rm Li}_5(x\alpha/y\beta)
+{\rm Li}_5(x\alpha y\eta)
+{\rm Li}_5(x\alpha\beta/\eta)
+{\rm Li}_5(x\xi y\beta)
+{\rm Li}_5(x\xi/y\eta)\nonumber\\&&{}
+{\rm Li}_5(x\xi\eta/\beta)
+{\rm Li}_5(\alpha y\beta/\xi)
+{\rm Li}_5(\alpha/\xi y\eta)
+{\rm Li}_5(\alpha\eta/\xi\beta)\nonumber\\&&{}-9\Bigl\{
 {\rm Li}_5(x y)
+{\rm Li}_5(x\beta)
+{\rm Li}_5(x\eta)
+{\rm Li}_5(x/y)
+{\rm Li}_5(x/\beta)
+{\rm Li}_5(x/\eta)\nonumber\\&&{}
+{\rm Li}_5(\alpha y)
+{\rm Li}_5(\alpha\beta)
+{\rm Li}_5(\alpha\eta)
+{\rm Li}_5(\alpha/y)
+{\rm Li}_5(\alpha/\beta)
+{\rm Li}_5(\alpha/\eta)\nonumber\\&&{}
+{\rm Li}_5(\xi y)
+{\rm Li}_5(\xi\beta)
+{\rm Li}_5(\xi\eta)
+{\rm Li}_5(y/\xi)
+{\rm Li}_5(\beta/\xi)
+{\rm Li}_5(\eta/\xi)\Bigr\}\nonumber\\&&{}+18\Bigl\{
 {\rm Li}_5(x)
+{\rm Li}_5(\alpha)
+{\rm Li}_5(\xi)
+{\rm Li}_5(y)
+{\rm Li}_5(\beta)
+{\rm Li}_5(\eta)
-{\rm Li}_5(1)\Bigr\}=\nonumber\\&&
\dfrac{3}{10}\log^5\xi+\dfrac34\left\{\log y-\log x\right\}\log^4\xi
+\dfrac32\left\{3\log y-\log\eta\right\}\log^2\eta\log^2\xi\nonumber\\&&{}
+\dfrac{1}{2}\pi^2\left\{\log\xi-3\log\eta\right\}\log^2\xi
+\dfrac{1}{5}\pi^4\log\xi\label{li5}
\end{eqnarray}
with $\xi:=1-x$, $\eta:=1-y$, $\alpha:=-x/\xi$, and $\beta:=-y/\eta$.
This formula is valid for all complex pairs $(x,y)$ in the neighbourhoods
of the 4 points of interest, namely
$(\frac12,\frac12)$, $(i,i)$, $(\frac12,i)$, and $(\frac12,-i)$.
Apart from the duplication and inversion relations~\cite{Lew}
\begin{equation}
{\rm Li}_5(-x)=-{\rm Li}_5(x)+\dfrac{1}{16}{\rm Li}_5(x^2)
={\rm Li}_5(-1/x)-\dfrac{1}{120}\log^5x-\dfrac{1}{36}\pi^2\log^3x
-\dfrac{7}{360}\pi^4\log x\label{ri}
\end{equation}
no further functional information about ${\rm Li}_5$ appears to be available.

Unfortunately, these functional identities do not deliver
all of the expected goods, since they yield only 4 independent
relations for the target set of polylogarithms.
One gives $\beta(5)=5\pi^5/1536$.
In terms of $\lambda(5)=\frac{31}{32}\zeta(5)$,
the other 3 give
\begin{eqnarray}
\widetilde{C}_5&=&
\dfrac{13}{81}\lambda(5)\label{r5c}\\
\widetilde{B}_5+\dfrac92\widetilde{D}_5&=&
\dfrac{47}{6}\lambda(5)\label{r51}\\
\widetilde{B}_5-\left(\dfrac{9}{2}\right)^3\widetilde{D}_5
+\left(\dfrac{5}{2}\right)^4\widetilde{E}_5&=&18\lambda(5)\label{r52}
\end{eqnarray}
the first of which was known~\cite{BBP,Lew}.
Resorting to numerical computation we found that
\begin{eqnarray}
\widetilde{B}_5&=&\dfrac{69}{8}\lambda(5)\label{qef}\\
\widetilde{H}_5&=&-\dfrac{1567}{3125}\beta(5)\label{n5h}
\end{eqnarray}
to an accuracy of 1,000 hexadecimal digits, which leaves no reasonable doubt
that these results are exact. Subsequently we found an intricate proof
of~(\ref{qef}), which will be given in section~3.
Hence it is proven that
\begin{equation}
\lambda(5)=\dfrac{31}{32}\zeta(5)=
\dfrac{8}{69}\widetilde{B}_5=
\dfrac{81}{13}\widetilde{C}_5=
-\dfrac{108}{19}\widetilde{D}_5=
-\dfrac{1250}{213}\widetilde{E}_5\label{r5}
\end{equation}
which can be solved to find SC$^*$ expressions for
all 4 of the constants $\zeta(5)$, $\log^52$, $\pi^2\log^32$,
and $\pi^4\log2$. In particular, we obtain
\begin{eqnarray}
\dfrac{2021}{256}\log^52&=&S_{5,1}(2783,-261592,-2783,-1500376,
-2783,-261592,2783,26717696)\nonumber\\
&&{}+S_{5,3}(29537,79446,-29537,-108983,-29537,79446,
29537,-49909)\nonumber\\
&&{}-26398\,S_{5,5}(1,0,-1,-1,-1,0,1,1)\label{l5}\\
\dfrac{62651}{2048}\zeta(5)&=&
9\,S_{5,1}(31,-1614,-31,-6212,-31,-1614,31,74552)\nonumber\\
&&{}+7\,S_{5,3}(173,284,-173,-457,-173,284,173,-111)\nonumber\\
&&{}-738\,S_{5,5}(1,0,-1,-1,-1,0,1,1)\,.\label{z5}
\end{eqnarray}

\subsection{Polylogarithms of order 6}

At $n=6$, there is a clear expectation that
$\widetilde{H}_6$ will combine with $\beta(5)\log2$,
in proportions predicted by~(\ref{n5h}), to yield a rational multiple of
$\beta(6)$.
Numerical computation gives
\begin{equation}
\frac{61\beta(6)}{3}=\frac{1567\beta(5)\log2
-3125\widetilde{H}_6}{2^8}\label{b6}
\end{equation}
which has been checked at high precision.
Clearly the three $\beta$--generating ladders~(\ref{fn}--\ref{hn})
are now exhausted.
Equally clearly, there is further mileage in the
$\lambda$--generating ladders~(\ref{an}--\ref{en}).

The general procedure is as follows.
Having found $7-k$ ladders that yield rational multiples of
$\lambda(2k-1)$ at order $n=2k-1$,
one forms $6-k$ combinations that vanish at this order.
These are expected to be rational multiples of $\pi^{2k}$ at order $n=2k$.
Having found these rational numbers by numerical computation,
one subtracts the same multiples of $\pi^{2k}L_{n-2k}$ at any subsequent
order, $n$. Thus one has $6-k$ ladders that vanish for $n\le2k$ and
are expected to give rational multiples of $\lambda(2k+1)$ at $n=2k+1$.
Having found these rational numbers by numerical computation,
one iterates the procedure, until it terminates with
a single rational multiple of $\lambda(11)=\frac{2047}{2048}\zeta(11)$.

Thus we should now form 3 ladders that vanish for $n\le6$, before progressing
to $n=7$. The following serve
\begin{eqnarray}
U_n&:=&\dfrac{13}{23}\widetilde{B}_n-\dfrac{243}{8}\widetilde{C}_n
-\dfrac{11041}{2048}\zeta(6)L_{n-6}\label{un}\\
V_n&:=&\dfrac{19}{23}\widetilde{B}_n+\dfrac{81}{2}\widetilde{D}_n
-\dfrac{87101}{12288}\zeta(6)L_{n-6}\label{vn}\\
W_n&:=&\dfrac{71}{23}\widetilde{B}_n+\dfrac{625}{4}\widetilde{E}_n
-\dfrac{1193757}{40960}\zeta(6)L_{n-6}\label{wn}
\end{eqnarray}
with combinations of the $n=5$ ladders
determined by~(\ref{r5}), and $\zeta(6)$ terms determined
by numerical computation at $n=6$. These constructs vanish for
$n\le6$.

\subsection{Polylogarithms of order 7}

Evaluating~(\ref{un}--\ref{wn}) at $n=7$, we find
rational multiples of $\lambda(7)=\frac{127}{128}\zeta(7)$, with
\begin{equation}
\dfrac{340}{23}\lambda(7)
=\dfrac{384}{463}U_7=\dfrac{32}{53}V_7=\dfrac{125}{819}W_7\label{z7}
\end{equation}
giving 3 independent SC$^*$ combinations of
$\{\zeta(7),\pi^6\log2,\pi^4\log^32,\pi^2\log^52,\log^72\}$.
These relations were
found at 64--bit precision and checked at much
greater precision.

\subsection{Polylogarithms of order 8}

Following the same systematic procedure, we find that
\begin{eqnarray}
X_n&:=&           463V_n-             636U_n-\dfrac{1323636287}{1769472}
\zeta(8)L_{n-8}\label{xn}\\
Y_n&:=&\dfrac{91}{25}V_n-\dfrac{265}{288}W_n-\dfrac{602893337}{113246208}
\zeta(8)L_{n-8}\label{yn}
\end{eqnarray}
vanish for $n\leq8$. The relative coefficients of previous ladders were
taken from~(\ref{z7}).

\subsection{Polylogarithms of order 9}

Evaluating~(\ref{xn},\ref{yn}) at $n=9$, we find
rational multiples of $\lambda(9)=\frac{511}{512}\zeta(9)$, with
\begin{equation}
\dfrac{217}{864}\lambda(9)
=\dfrac{1}{10435}X_9
=\dfrac{500}{37403}Y_9
\label{z9}
\end{equation}
giving two independent SC$^*$ combinations of $\zeta(9)$ and products of
powers of $\pi$ and $\log2$.

\subsection{Polylogarithms of order 10}

Proceeding as before, we find that
\begin{equation}
Z_n:=\dfrac{1}{4823}\left(2087Y_n-\dfrac{37403}{2500}X_n\right)
-\dfrac{12227440999}{135895449600}\zeta(10)L_{n-10}\label{zn}
\end{equation}
vanishes for $n\le10$. The relative coefficients of previous ladders were
taken from~(\ref{z9}).

\subsection{Polylogarithms of order 11}

The process terminates with the numerical result that
\begin{equation}
\lambda(11)=\dfrac{2047}{2048}\zeta(11)
=\dfrac{129600000}{41323873}Z_{11}\label{z11}
\end{equation}
enabling one to evaluate in SC$^*$ a combination of $\zeta(11)$
with products of powers of $\pi$ and $\log2$. At no stage
of the iterative process does one need more than 128--bit precision
to be very confident of the single rational number that is
found at each step of each surviving ladder.  Yet, unwrapping~(\ref{z11}),
one obtains the arcane integer relation
\begin{eqnarray}&&
{\tt 46090055410032553920}\,\zeta(11)=
{\tt 105497707483968307200}\,\Re\,{\rm Li}_{11}(\dfrac{1+i}{2})\nonumber\\
&&{}+{\tt 14102390469191270400}\,\Re\,{\rm Li}_{11}(\dfrac{1+i}{4})
-{\tt 943412955347681280}\,\Re\,{\rm Li}_{11}(\dfrac{1+i}{8})\nonumber\\&&{}
+{\tt 8628616191131674214400}\,{\rm Li}_{11}(\dfrac12)
+{\tt 8666542920405771878400}\,{\rm Li}_{11}(-\dfrac12)\nonumber\\&&{}
+{\tt 8389140238437235200}\,{\rm Li}_{11}(-\dfrac14)
-{\tt 73384332676300800}\,{\rm Li}_{11}(-\dfrac18)\nonumber\\&&{}
-{\tt 5097106123776}\log^{11}2
+{\tt 9394465639680}\,\pi^2\log^{9}2\nonumber\\&&{}
-{\tt 13065007342464}\,\pi^4\log^{7}2
+{\tt 20585306545056}\,\pi^6\log^{5}2\nonumber\\&&{}
-{\tt 42801564610332}\,\pi^8\log^{3}2
+{\tt 139087141363625}\,\pi^{10}\log2\label{f11}
\end{eqnarray}
which has been checked at 16000--bit precision. Working at this
precision, the lattice algorithm {\sc pslq}~\cite{PSLQ} proved
that between the 13 constants on the right there exists
no relation with integer coefficients of less than 300 decimal digits.
Thus one may be confident that~(\ref{f11}) is both correct and unique.
Hence we reach the terminus suggested by~(\ref{w},\ref{h}).

\section{Hypergeometric series and Euler sums}

It was noted in section~2.5 that~(\ref{r5c}--\ref{r52}) were proven
by using~(\ref{li5}), but~(\ref{qef}) was not. A proof
is given here; it illustrates an important connection, via
hypergeometric series~\cite{WP2},
to Euler sums~\cite{EUL,BBG,BG,DJB,BBB,MEH,FS}.

First we use the integral representation~\cite{Lew}
\begin{equation}
{\rm Li}_n(z)=\int_0^1\frac{dx}{z^{-1}-1+x}\,
\frac{(-\log(1-x))^{n-1}}{(n-1)!}
\label{irep}
\end{equation}
to obtain generating functions for~(\ref{an}--\ref{cn})
\begin{eqnarray}
\sum_{n>0}A_nt^n&=&t\int_0^1\frac{dx}{(1-x)^t}
\,\frac{d\log(1+x)}{dx}\label{ag}\\
2\sum_{n>0}B_nt^n&=&t\int_0^1\frac{dx}{(1-x)^{2t}}
\,\frac{d\log(1+x^2)}{dx}\label{bg}\\
2\sum_{n>0}C_nt^n&=&t\int_0^1\frac{dx}{(1-x)^t}
\,\frac{d\log(\frac{3+x^2}{3-x})}{dx}\label{cg}\\
2\sum_{n>0}D_nt^n&=&t\int_0^1\frac{dx}{(1-x)^{2t}}
\,\frac{d\log(\frac{(1+x)^2+x^2(1-x)^2}{1+x^2})}{dx}\label{dg}
\end{eqnarray}
that make the $n=1$ results, $A_1=2B_1=2C_1=2D_1=\log2$, immediately
apparent.

\subsection{A hypergeometric proof that $\zeta(5)$ is in SC$^*$}

To reduce the $p=1$ cases~(\ref{ag},\ref{bg}) to ${}_3F_2$ series,
we change variables
to $u:=4x/(1+x)^2$, in the former, and to $v:=2x/(1+x^2)$, in the latter,
obtaining
\begin{eqnarray}
\sum_{n>0}A_nt^n&=&\frac{t}{4}\int_0^1\frac{du}{(1-u)^{t/2}}\,f(u,t)
\label{agu}\\
2\sum_{n>0}B_nt^n&=&\frac{t}{2}\int_0^1\frac{v\,dv}{(1-v)^{t}}\,f(v^2,t)
\label{bgv}\\
f(z,t)&:=&(1-z)^{-1/2}\left(\dfrac12+\dfrac12(1-z)^{1/2}\right)^{t-1}
\nonumber\\
&=&\,{}_2F_1(\dfrac32-\dfrac12t,\,1-\dfrac12t;\,2-t;\,z)\label{fzt}
\end{eqnarray}
with hypergeometric integrands given by Eq~(15.1.13) of~\cite{AS}.
Expanding~(\ref{fzt}) and integrating~(\ref{agu},\ref{bgv}), we obtain
\begin{eqnarray}
\sum_{n>0}A_nt^n&=&1-{}_3F_2(-\dfrac12t,\,\dfrac12-\dfrac12t,\,1;\,
1-\dfrac12t,\,1-t;\,1)\label{ais}\\
2\sum_{n>0}B_nt^n&=&1-{}_3F_2(-\dfrac12t,\,\dfrac12,\,1;\,
1-\dfrac12t,\,1-t;\,1)\label{bis}
\end{eqnarray}
where~(\ref{ais}) was obtained by similar means in~\cite{WP2}.
Thus, to prove~(\ref{qef}),
we need to derive expansions to $O(t^5)$ of series of the form
\begin{equation}
F(a,b,c):=1-{}_3F_2(-a,\,\dfrac12-b,\,1;\,
1-a,\,1-c;\,1)\label{tpc}
\end{equation}
where $a,b,c$ are rational multiples of $t$ specified
by~(\ref{ais},\ref{bis}).

Such hypergeometric series were intensively investigated
in~\cite{WP2}, in connection
with integrals arising from the Feynman diagrams of perturbative
quantum field theory. It was shown that
there is a wreath product group,
$S_3\wr Z_2$, that relates 72 Saalsch\"utzian series of the form
${}_3F_2(-\alpha_1,-\alpha_2,1;1+\alpha_3,1+\alpha_4;1)$
via linear transformations of the 4 parameters.
In particular,~(\ref{tpc}) and related series,
of importance in quantum field theory~\cite{DJB}
and its relation~\cite{BKP,BGK,BK15} to knot theory,
may be transformed to series of the more symmetrical type~\cite{WP2}
\begin{equation}
W(a_1,a_2;a_3,a_4):=(\dfrac12+a_3)^{-1}(\dfrac12+a_4)^{-1}\,
{}_3F_2(\dfrac12-a_1,\,\dfrac12-a_2,\,1;\,
\dfrac32+a_3,\,\dfrac32+a_4;\,1)\label{ww}
\end{equation}
with $W(0,0;0,0)=4\lambda(2)=\frac12\pi^2$. The group of transformations
that constrains the expansion of~(\ref{ww}) in small parameters
$a_k$ is the symmetry group of the square, resulting
from the obvious symmetries in $(a_1,a_2)$ and $(a_3,a_4)$ and the
non-trivial reduction~\cite{WP2}
\begin{equation}
W(a_1,a_2;a_3,a_4)+W(a_3,a_4;a_1,a_2)
=\frac{\Gamma\left(1+\sum_{k}a_k\right)}{\prod_k\Gamma(\frac12+a_k)}
\prod_{i=1,2}\prod_{j=3,4}B(\dfrac12+a_i,\dfrac12+a_j)
\label{inv}
\end{equation}
with $B(a,b):=\Gamma(a)\Gamma(b)/\Gamma(a+b)$.
The transformation from~(\ref{tpc}) to~(\ref{ww}) is
\begin{equation}
F(a,b,c)=1-\frac{\pi a\,B(\frac12+b,1-c)}
{\sin(\pi a)\,B(\frac12-a+b,1+a-c)}
+ a c W(c-b,a-b;b,b-c)\,.\label{fab}
\end{equation}
In particular,~(\ref{ais},\ref{bis}) transform to
\begin{eqnarray}
\sum_{n>0}A_nt^n=F(\dfrac12t,\dfrac12t,t)
&=&1-\frac{\pi t}{2^t\sin(\pi t)}
+\dfrac12t^2W(\dfrac12t,0;\dfrac12t,-\dfrac12t)\label{cona}\\
&=&1-\frac{\pi t}{2^t\tan(\pi t)}
-\dfrac12t^2W(\dfrac12t,-\dfrac12t;\dfrac12t,0)\label{alta}\\
2\sum_{n>0}B_nt^n=F(\dfrac12t,0,t)
&=&1-\frac{\frac12\pi t}{2^t\sin(\frac12\pi t)}
+\dfrac12t^2W(t,\dfrac12t;0,-t)\label{conb}\\
&=&1-\frac{2\pi t\cos(\frac32\pi t)}{2^t\sin(2\pi t)}
-\dfrac12t^2W(0,-t;t,\dfrac12t)\label{altb}
\end{eqnarray}
with the beta functions of~(\ref{fab}) yielding a factor of $2^{-t}$ and
trigonometric functions in~(\ref{cona},\ref{conb}) and in the
alternative forms~(\ref{alta},\ref{altb}) provided by~(\ref{inv}).

The task of expanding~(\ref{conb}) to $O(t^5)$
is not so severe as might appear, thanks to the $t^2$ in front
of $W$, whose expansion is highly constrained by symmetries.
From the analysis of~\cite{WP2} we obtain the general expansion
\begin{eqnarray}
\dfrac18t^2W(a_1t,a_2t;a_3t,a_4t)-\dfrac18t^2W(a_3t,a_4t;a_1t,a_2t)=
\delta_1\sum_{n\ge3}(2\sigma_1)^{n-3}\overline{A}_nt^n\nonumber\\
\quad{}+\left(\delta_2-\dfrac{49}{45}\sigma_1\delta_1\right)
2^{-5-2\sigma_1 t}(\pi t)^4+f_5(a_1,a_2,a_3,a_4)\lambda(5)t^5
+O(t^6)\label{asy}
\end{eqnarray}
with $\sigma_n:=\sum_k a_k^n$, $\delta_1:=a_1+a_2-a_3-a_4$,
$\delta_2:=a_1a_2-a_3a_4$ and
\begin{equation}
f_5(a_1,a_2,a_3,a_4):=2\sigma_2\delta_1
-3\sigma_1(\delta_2+\sigma_1\delta_1)\label{f5}
\end{equation}
giving $f(\frac12,0,\frac12,-\frac12)=0$, as required by~(\ref{cona}).
Using~(\ref{asy}) in conjunction with the elementary expansion
of~(\ref{inv}), one may expand~(\ref{fab}) to $O(t^5)$.
In particular,~(\ref{conb}) gives
\begin{equation}
\widetilde{B}_5=f_5(1,\dfrac12,0,-1)\lambda(5)
=\dfrac{69}{8}\lambda(5)\label{bl5}
\end{equation}
which finally proves~(\ref{qef}) and hence~(\ref{r5}--\ref{z5}).

\subsection{Euler sums from hypergeometric series}

We note that~(\ref{cona},\ref{conb})
establish the reducibility of~(\ref{an},\ref{bn}), at any order $n$,
to alternating Euler sums~\cite{DJB}.
The results of sections~2.1--2.5 then establish
that~(\ref{cn}--\ref{en}) are reducible to Euler sums for $n\le5$.
Moreover, {\sc pslq} found such reductions
at $n=6,7$ and, more generally, indicates that
a ladder giving a multiple of $\zeta(2k+1)$ at order $n=2k+1$
gives Euler sums of depth $n-2k$ at order $n>2k+1$. This suggests
that the generators~(\ref{cg},\ref{dg}) might also be of the
form~(\ref{tpc}). We thus worked backwards, from the observations that
\begin{eqnarray}
\widetilde{C}_5=\dfrac{13}{81}\lambda(5)
=\dfrac23f_5(\dfrac12,\dfrac13,\dfrac16,-\dfrac12)\lambda(5)\label{cl5}\\
\widetilde{D}_5=-\dfrac{19}{108}\lambda(5)
=\dfrac23f_5(\dfrac13,\dfrac16,\dfrac13,-\dfrac13)\lambda(5)\label{dl5}
\end{eqnarray}
to prove the remarkably simple results
\begin{eqnarray}
F(\dfrac12t,\dfrac16t,\dfrac23t)&=&2\sum_{n>0}C_nt^n+O(t^6)\label{remc}\\
F(\dfrac12t,\dfrac13t,\dfrac23t)&=&2\sum_{n>0}D_nt^n+O(t^6)\label{remd}
\end{eqnarray}
which strongly suggest that these are full generators.
We thus propose that
\begin{eqnarray}
2\sum_{n>0}C_nt^n&=&1-
{}_3F_2(-\dfrac12t,\,\dfrac12-\dfrac16t,\,1;\,1-\dfrac12t,\,1-\dfrac23t;\,1)
\label{cis}\\
2\sum_{n>0}D_nt^n&=&1-
{}_3F_2(-\dfrac12t,\,\dfrac12-\dfrac13t,\,1;\,1-\dfrac12t,\,1-\dfrac23t;\,1)
\,.\label{dis}
\end{eqnarray}

To prove that~(\ref{cis},\ref{dis}) are true,
one may compare the poles of the left and right--hand sides,
which occur only on the positive real axis.
On the left, we have
\begin{eqnarray}
2\sum_{n>0}C_nt^n
&=&t\sum_{k>0}(-1)^k\left\{\frac{2^{-3k}}{k-t/3}
-\frac{2^{1-k}}{k-t}\right\}\label{leftc}\\
2\sum_{n>0}D_nt^n
&=&t\sum_{k>0}\left\{\frac{2^{1-3k/2}\cos(\frac14\pi k )}{k-2t/3}
-\frac{(-\frac14)^{k}}{k-t/2}\right\}\label{leftd}
\end{eqnarray}
by virtue of the definitions~(\ref{cn},\ref{dn}). On the right, we have
\begin{eqnarray}
F(\dfrac12t,\dfrac16t,\dfrac23t)
&=&1-\frac{\frac12\pi t}{2^t\sin(\frac12\pi t)\cos(\frac16\pi t)}
+\dfrac13t^2W(\dfrac12t,\dfrac13t;\dfrac16t,-\dfrac12t)\label{conc}\\
&=&1-\frac{\pi t}{2^t\tan(\pi t)\cos(\frac13\pi t)}
-\dfrac13t^2W(\dfrac16t,-\dfrac12t;\dfrac12t,\dfrac13t)\label{altc}\\
F(\dfrac12t,\dfrac13t,\dfrac23t)
&=&1-\frac{\frac12\pi t}{2^t\sin(\frac12\pi t)\cos(\frac13\pi t)}
+\dfrac13t^2W(\dfrac13t,\dfrac16t;\dfrac13t,-\dfrac13t)\label{cond}\\
&=&1-\frac{\frac12\pi t\cos(\frac56\pi t)}
{2^t\sin(\frac12\pi t)\cos(\frac16\pi t)\cos(\frac13\pi t)}
-\dfrac13t^2W(\dfrac13t,-\dfrac13t;\dfrac13t,\dfrac16t)\label{altd}
\end{eqnarray}
with~(\ref{conc},\ref{cond}) obtained from~(\ref{fab}), and the alternative
forms~(\ref{altc},\ref{altd}) from~(\ref{inv}).
Elementary analysis of the trigonometric parts of~(\ref{altc},\ref{altd})
reveals that their poles, on the positive real axis,
coincide with those of~(\ref{leftc},\ref{leftd}). Since the $W$ series
of~(\ref{altc},\ref{altd}) are finite on the positive real axis,
we conclude that left and right--hand sides of~(\ref{cis},\ref{dis})
differ, if at all, by entire functions.
To prove that they are equal, it suffices to show that they give
the same values at infinity. On the left,~(\ref{leftc},\ref{leftd})
give pairs of geometric series at infinity;
on the right, the hypergeometric series become geometric. Hence one has
merely to verify that
\begin{eqnarray}
-3\sum_{k>0}\,(-\dfrac18)^k+2\sum_{k>0}\,(-\dfrac12)^k=
\dfrac13-\dfrac23=-\dfrac13=1-\sum_{k\ge0}\,(\dfrac14)^k\label{geoc}\\
-3\Re\sum_{k>0}\,(\dfrac{1+i}{4})^k+2\sum_{k>0}\,(-\dfrac14)^k=
-\dfrac35-\dfrac25=-1=1-\sum_{k\ge0}\,(\dfrac12)^k\label{geod}
\end{eqnarray}
to conclude that the differences between left and right
are entire functions that vanish at infinity and hence at all $t$.

These proofs give little clue as to the origin
of the remarkable results~(\ref{cis},\ref{dis}); they
merely certify what was already strongly suggested by
observations~(\ref{remc},\ref{remd}). If one had a method of systematic
derivation, as opposed to mere proof, the outstanding problem of the
$p=5$ case~(\ref{en}) might be more tractable.

\subsection{Catalan's constant from hypergeometric series}

From series~(\ref{bis},\ref{dis}) we can derive hypergeometric
generators of $F_n$ and $G_n$, yielding $G:=\beta(2)$ at $n=2$.
Defining
\begin{eqnarray}
F(t):=2\sum_{n>0}(B_n+iF_n)t^n&=&
2t\sum_{k>0}\frac{\left(\frac{1+i}{2}\right)^k}{k-2t}\label{bif}\\
G(t):=2\sum_{n>0}(D_n+iG_n)t^n&=&
2t\sum_{k>0}\left\{\frac{\left(\frac{1+i}{4}\right)^k}{k-2t/3}
-\frac{\left(\frac{-i}{2}\right)^k}{k-t}\right\}\label{dig}\\
H(t):=2\sum_{n>0}(E_n+iH_n)t^n&=&
2t\sum_{k>0}\left\{\frac{\left(\frac{1-i}{8}\right)^k}{k-2t/5}
-\frac{2\left(\frac{-i}{2}\right)^k}{k-t}\right\}\label{eih}
\end{eqnarray}
we obtain the recurrence relations
\begin{eqnarray}
\frac{2\,F(t)}{t}-\frac{i\,F(t-1)-i}{t-1}&=&
\frac{2+2i}{1-2t}\label{frec}\\
\frac{2^3G(t)}{t}-\frac{i\,G(t-3)-i}{t-3}&=&
\frac{12+12i}{3-2t}+\frac{8i}{1-t}+\frac{4}{2-t}\label{grec}\\
\frac{2^5H(t)}{t}+\frac{i\,H(t-5)-i}{t-5}&=&
\frac{40-40i}{5-2t}+\frac{64i}{1-t}+\frac{32}{2-t}
-\frac{16i}{3-t}-\frac{8}{4-t}\label{hrec}
\end{eqnarray}
and can hence relate the generators of $F_n$, $G_n$ and $H_n$
to those of $B_n$, $D_n$ and $E_n$,
by taking the imaginary parts of~(\ref{frec}--\ref{hrec}).
In particular the proven results~(\ref{bis},\ref{dis}) yield
\begin{eqnarray}
\Im F(t):=2\sum_{n>0}F_nt^n&=&t(1-t)^{-1}\,
{}_3F_2(\dfrac12-\dfrac12t,\,\dfrac12,\,1;\,
\dfrac32-\dfrac12t,\,1-t;\,1)\label{fis}\\
\Im G(t):=2\sum_{n>0}G_nt^n&=&t(1-t)^{-1}\,
{}_3F_2(\dfrac12-\dfrac12t,\,\dfrac12-\dfrac13t,\,1;\,
\dfrac32-\dfrac12t,\,1-\dfrac23t;\,1)\label{gis}
\end{eqnarray}
with the $O(t)$ terms giving $2F_1=2G_1={}_2F_1(\frac12,
\frac12;\frac32;1)=\frac12\pi=2\beta(1)$.

It was shown in section~2.2 that Catalan's constant
\begin{equation}
G=\dfrac32(F_2-G_2)=\Im\left\{3\,{\rm Li}_2(\dfrac{1+i}{2})
-{\rm Li}_2(\dfrac{1+i}{4})+\dfrac32\,{\rm Li}_2(\dfrac{-i}{2})\right\}
\label{cat}
\end{equation}
is obtained at $n=2$, with~(\ref{G}) giving the corresponding SC$^*$ series.
As far as we can tell,~(\ref{cat}) is a new result;
we were unable to locate it in Victor Adamchik's interesting
compilation~\cite{AA} of representations for Catalan's constant.
By expanding~(\ref{fis},\ref{gis}) to $O(t^2)$,
we may transform it to
\begin{equation}
G=\sum_{n=1}^\infty\frac{\left(\frac12\right)^{2n+1}}{2n+1}
{2n\choose n}\sum_{m=1}^{2n}\frac{1}{m}
\label{acat}
\end{equation}
which is vastly inferior to~(\ref{cat}), for computational purposes,
yet serves to illustrate the type of non-eulerian
double sum that results from the shifts in~(\ref{frec},\ref{grec}).

It is straightforward to determine whether
${}_3F_2(a_1+b_1t,a_2+b_2t,1;a_3+b_3t,a_4+b_4t;1)$
has a small-$t$ expansion yielding
non-alternating Euler sums, alternating Euler sums, or non-eulerian sums,
for values of $a_k$ are integers (denoted by $n$) or half--integers
(denoted by $h$). The sole pattern $(a_1,a_2)(a_3,a_4)$ that yields
non-alternating Euler sums is $(n,n)(n,n)$.
Alternating Euler sums result from the patterns $(h,n)(n,n)$, $(n,n)(h,n)$,
$(h,n)(h,n)$ and $(h,h)(h,h)$. Non-eulerian sums like~(\ref{acat})
result from the remaining patterns: $(h,n)(h,h)$, $(h,h)(h,n)$,
$(h,h)(n,n)$ and $(n,n)(h,h)$. Thus the shifts
in~(\ref{frec},\ref{grec}) transform the generators~(\ref{bis},\ref{dis})
of the alternating Euler-sums in $B_n$ and $D_n$ to the
generators~(\ref{fis},\ref{gis}) of the non-eulerian sums in $F_n$ and $G_n$.
The shift in~(\ref{hrec}) presumably accomplishes the same transformation
from the Euler sums in $E_n$ to the non-eulerian sums in $H_n$, of which
$\beta(5)$ in~(\ref{n5h}) is the sole discovery at $n\le5$ that remains
unproven. Hence the outstanding challenge is to find a hypergeometric
generator for the Euler sums in $E_n$, which can then be transformed,
via~(\ref{hrec}), to prove results such as~(\ref{n5h},\ref{b6}) in the
non-eulerian sector.

\subsection{Data on the remaining case}

It is proven, in sections~3.1--3.2, that~(\ref{an}--\ref{dn}) are reducible
to Euler sums, at any order $n$, since the $O(t^n)$ term in~(\ref{ww})
yields only Euler sums. From the discovery~(\ref{z11}), we then infer that
$E_n$ of~(\ref{en}) is reducible to Euler sums for $n\le11$.
Moreover, {\sc pslq} found that~(\ref{zn}) is reducible to
alternating double sums at $n=12$, which implies that $E_{12}$
reduces to Euler sums of depths up to 10. Assuming that nothing
untoward happens at $n\ge13$, we infer that there is a hypergeometric
generator for $E_n$ akin to, but probably more complicated than,
those found for $A_n$ and $B_n$ in~(\ref{cona}--\ref{altb})
and for $C_n$ and $D_n$ in~(\ref{conc}--\ref{altd}).

Accordingly, a hypergeometric representation of
\begin{equation}
E(t):=\Re H(t)=2\sum_{n>0}E_nt^n=
t\sum_{k>0}\left\{\frac{2^{1-5k/2}\cos(\frac14\pi k)}{k-2t/5}
-\frac{2(-\frac14)^k}{k-t/2}\right\}\label{what}
\end{equation}
was earnestly sought. We may write it as
\begin{equation}
E(t)=1-\frac{\frac12\pi t}{2^t\sin(\frac12\pi t)}
\left(\frac{1}{\cos(\dfrac15\pi t)}-8\sin^2(\dfrac15\pi t)\right)
-\dfrac25U(t)\label{U}
\end{equation}
with a trigonometric part that removes poles
on the positive real axis, since
\begin{eqnarray}
\frac{1}{\cos(\dfrac15\pi t)}-8\sin^2(\dfrac15\pi t)
&=&\left\{\begin{array}{rl}
 1&\mbox{for $t=0$~mod~10}\\
-4&\mbox{for $t=2,4,6,8$~mod~10}
\end{array}\right.\label{emod}\\
\frac{1}{\sin(\frac12\pi t)}&=&
-2\sin(\dfrac1{10}\pi t)\quad\mbox{for $2t=5$~mod~10}\label{omod}
\end{eqnarray}
provide the correct residues. Moreover, the subtraction made
in~(\ref{U}) is the unique trigonometric term that both
removes poles and also vanishes at infinity in the right half--plane,
as was the case for the corresponding terms in~(\ref{altc},\ref{altd}).

By construction, $U(t)$ in~(\ref{U}) is finite for $t>-2$.
It is positive for $t>0$, increasing from $U(0)=0$, through
$U(5)=\frac{20}{3}$, to a maximum value $U(t_{\rm max})\approx6.786$ at
$t_{\rm max}\approx 6.731$, and then falling, through $U(10)=\frac{20}{3}$,
to $U(\infty)=6$.
A great deal of further data on $U(t)$ is available.
Its poles on the negative real axis are determined by the fact that $E(t)$
is finite there. Its expansion around $t=0$ involves only Euler sums
to $O(t^{12})$. The results to order $n=6$ are
\begin{eqnarray}
U(t)&=&\frac{(\pi t)^2+\frac{53}{1200}(\pi t)^4
+\frac{30131}{25200000}(\pi t)^6}{2^{1+t}}
-\frac65\sum_{n\ge3}\overline{A}_nt^n\nonumber\\&&{}
+\frac{213}{250}\left\{\frac{31}{32}\zeta(5)t^5+\left(\frac{3}{16}\zeta^2(3)
-\frac12\sum_{m>n>0}\frac{(-1)^{m+n}}{m^5n}\right)t^6\right\}+O(t^7)
\label{expu}
\end{eqnarray}
where $\overline{A}_n$ absorbs Euler sums of depth $n-2$ at order $n$,
leaving sums of depths up to $n-4$, with the first alternating double sum
appearing at $n=6$.

From the recurrence relation~(\ref{hrec}),
we deduce that $U(t)$ has a rational asymptotic expansion, beginning with
\begin{equation}
U(t)\sim6\left\{1+\dfrac{11}{10t}+\dfrac{157}{(10t)^2}
-\dfrac{1749}{(10t)^3}-\dfrac{433651}{(10t)^4}
-\dfrac{43430405}{(10t)^5}-\dfrac{4000517955}{(10t)^6}+O(t^{-7})\right\}\,.
\label{uas}
\end{equation}
Defining $k_n$ as the coefficient of $6/(10t)^{n}$, we find
that for $1\le n\le1000$ it is an odd integer, divisible by 3 if $n=0$~mod~3,
by 5 precisely $\sum_{k>0}\lfloor5^{-k}n\rfloor$ times,
by 7 if $n=0$~mod~6, by 11 if $n=0,1,3$~mod~10, and by $n+1$ if $n+1$
is a prime greater than 11. No further pattern of factors is apparent.
For example, the 84-digit prime
\begin{equation}
P_{84}=\mbox{\footnotesize$
20464734789428471449753650289585781678688420563295
0514231255723964171455482439213639$}\label{p84}
\end{equation}
remains in $k_{39}$.
The sign of $k_n$ changes at $n=3,11,18,25,33,40\ldots$, with runs of
7 or 8 coefficients of the same sign, giving an average interval
$\overline\Delta\approx7.38$.
These oscillations may be understood by an application of Cauchy's theorem,
which relates the contour integral
\begin{equation}
\frac54\int_{-\infty}^{\infty}\frac{dx\,x\exp(-i x\log2)}{(z-ix)
\sinh(\frac12\pi x)}\left({1\over\cosh(\frac15\pi x)}
+8\sinh^2(\dfrac15\pi x)\right)=\left\{\begin{array}{lr}
12-2U(z)&\mbox{for $\Re\,z>0$}\\
7+5E(z)&\mbox{for $\Re\,z<0$}
\end{array}\right.\label{cont}
\end{equation}
to $U(z)$ in the right half-plane and $E(z)$ in the left half-plane.
Thus the coefficients in~(\ref{uas}) are moments of the
subtraction term along the imaginary axis. For example,~(\ref{p84})
is
\begin{equation}
P_{84}=-{2^{36}\,5^{32}\over9}\int_0^\infty dx\,{x^{39}
\cos(x\log2)\over\sinh(\frac12\pi x)}\,\left({1\over\cosh(\frac15\pi x)}
+8\sinh^2(\dfrac15\pi x)\right)\label{m39}
\end{equation}
whose value was checked numerically, by expanding the hyperbolic terms in
powers of $\exp(-\frac1{10}\pi x)$. Similarly, $k_{906}$
contains the 3139--digit prime
\begin{equation}
P_{3139}={2^{903}\,5^{682}\over514269}\int_0^\infty dx\,{x^{906}
\sin(x\log2)\over\sinh(\frac12\pi x)}\,\left({1\over\cosh(\frac15\pi x)}
+8\sinh^2(\dfrac15\pi x)\right)\,.\label{m257}
\end{equation}
The average interval between sign changes in the asymptotic expansion
of~(\ref{cont}) is
\begin{equation}
\overline\Delta=\frac{\pi}{\arctan({\pi/\log1024})}
\approx7.38257\label{dbn}
\end{equation}
in good agreement with the changes observed up to $n=1000$.

We also obtain rational results when $t=0$~mod~5. On the positive
axis, one has
\begin{equation}
U(5n)=25n\,\Re\left\{
\sum_{k=1}^{2n}\frac{(2+2i)^k-2}{k}\left(\frac{i}{2}\right)^{5n-k}
-\sum_{k=2n+1}^{5n}\frac{2}{k}\left(\frac{i}{2}\right)^{5n-k}\right\}
\label{ratp}
\end{equation}
giving $U(5)=U(10)=20/3$. On the negative axis, the values at odd $n$ are
\begin{equation}
U(-5n)=V(n):=
-25n\,\Re\left\{
\sum_{k=1}^{2n-1}\frac{(2+2i)^{-k}-2}{k}\left(\frac{i}{2}\right)^{k-5n}
-\sum_{k=2n}^{5n-1}\frac{2}{k}\left(\frac{i}{2}\right)^{k-5n}\right\}
\label{ratm}
\end{equation}
giving $U(-5)=1900/3$.
At negative $t=0$~mod~10, where $U(t)$ is singular,
we nonetheless have a rational remainder, after subtraction of the
rational pole term, with
\begin{equation}
U(\varepsilon-10n)-(-2^{10})^n\left\{\frac{25n}{\varepsilon}-\frac52\right\}
=V(2n)+O(\varepsilon)\label{rats}
\end{equation}
giving $U(\varepsilon-10)=-25600/\varepsilon+20310+O(\varepsilon)$.

Finally, we have rational results for
\begin{equation}
\widetilde{U}(t):=U(t)-\frac{5\pi t}{2^t\sin(\dfrac12\pi t)}\label{ut}
\end{equation}
when $2t=5$~mod~10. For positive $t=5n/2$, with $n$ odd,
\begin{equation}
\widetilde{U}(5n/2)
=25n\sum_{k=0}^{n-1}\frac{\Re(4+4i)^{-k}}{2n-2k}
-50n\sum_{k=0}^{\lfloor\frac54n\rfloor}\frac{(-4)^{-k}}{5n-4k}\label{rathp}
\end{equation}
gives $\widetilde{U}(5/2)=15$. Near $t=-5n/2$, we have
a rational residue and remainder in
\begin{eqnarray}
&&\widetilde{U}(\varepsilon-5n/2)-\Re(4+4i)^n
\left\{\frac{125n}{4\varepsilon}-\frac{25}{2}\right\}=\nonumber\\&&
-25n\sum_{k=1}^{n-1}\frac{\Re(4+4i)^k}{2n-2k}
+50n\sum_{k=1}^{\lfloor\frac54n\rfloor}\frac{(-4)^k}{5n-4k}
+O(\varepsilon)\label{rathm}
\end{eqnarray}
with $n$ odd and positive,
giving $\widetilde{U}(\varepsilon-5/2)=125(1/\varepsilon-2)+O(\varepsilon).$

It is remarkable that simple trigonometric subtraction in~(\ref{U})
removes both $\pi$ and $\log 2$ from the results
of~(\ref{ratp}--\ref{rats}) for $t=0$~mod~5,
and that~(\ref{ut}) similarly gives rational results
in~(\ref{rathp},\ref{rathm}) for $2t=5$~mod~10.
Also of note is the factor $(2+2i)^k-2$ in the first term of~(\ref{ratp}).
The circumstance that $\Re(2+2i)^p=2$~mod~$p$, for all prime $p$, removes
primes from the denominator of $U(5n)$. For example, when $n$ is odd
and positive, no prime $p\in[5n/3,2n]$ can occur in the denominator, though
larger primes, $p\in[2n,5n]$, result from the second term of~(\ref{ratp}).
Such behaviour is highly specific to the combination~(\ref{what})
and might be expected if it involved a hypergeometric series that
terminates when $t=5$~mod~10. Indeed, the complex of results in
this section suggests that $U(t)$ may be expressible as a series with
parameters that are half-integer at $t=0$, degenerating to a terminating
series when $t/10$ is half-integer.

\subsection{Impasse}

Despite intensive investigation, we were unable to find a ${}_3F_2$
series of type~(\ref{ww}) that reproduces the wealth of data on $U(t)$.
The {\em proxime accessit\/} to the expansion~(\ref{expu}) is
\begin{equation}
U(t)=t^2W(\dfrac35t,-\dfrac12t;0,\dfrac25t)+
\dfrac{1}{50}2^{-t}\pi^2\zeta(3)t^5
-\left\{\dfrac{21}{1000}\zeta^2(3)+\dfrac{19}{70000}\pi^6\right\}t^6
+c_7t^7+O(t^8)\label{prox}
\end{equation}
in which the $W$ series, rather remarkably, leaves only products of
$\pi^2$, $\log 2$
and $\zeta(3)$ to $O(t^6)$. However, {\sc pslq} found that
an irreducible alternating triple Euler sum enters at $O(t^7)$,
in the coefficient
\begin{eqnarray}
c_7&=&
-\dfrac{99791}{120000}\zeta(7)
+\dfrac{71}{180000}\pi^6\log2
-\dfrac{31}{500}\zeta(5)\log^22
+\dfrac{5639}{60000}\zeta(5)\pi^2
+\dfrac{199}{180000}\zeta(3)\pi^4\nonumber\\&&{}
+\frac{2}{375}\sum_{k>m>n>0}\left\{6\frac{(-1)^{m+n}}{k^5m\,n}
-\frac{(-1)^{k+n}}{k^3m^3n}\right\}
\label{c7}
\end{eqnarray}
thereby dashing the hope that the corrections in~(\ref{prox}) might
have been generated by gamma functions, or their derivatives.
It thus appears that the simplicity of~(\ref{prox}) to $O(t^6)$,
which far outstrips any rival Ansatz involving a single $W$ series,
is illusory.

To appreciate the magnitude of the problem that we face, it is instructive
to consider the identities that underwrite the validity
of~(\ref{cona},\ref{conb},\ref{conc},\ref{cond}), which involve $W$
series that are singular on the positive real axis. These singularities
must be generated trigonometrically, since the
singularities of the generators are trigonometric. Consistency is ensured
by the peculiar identities
\begin{eqnarray}
\left.\frac{(\frac12)_n}{(\frac12+\frac12t)_n}
\right|_{n=\frac12t-\frac12}&=&2^{1-t}\label{poca}\\
\left.\frac{(\frac12-\frac12t)_n}{(\frac12)_n}
\right|_{n=t-\frac12}\,\,\,
&=&2^{1-t}\cos(\dfrac12\pi t)\label{pocb}\\
\left.\frac{(\frac12-\frac13t)_n}{(\frac12+\frac16t)_n}
\right|_{n=\frac12t-\frac12}&=&2^{2-t}\cos(\dfrac13\pi t)\label{pocc}\\
\left.\frac{(\frac12-\frac16t)_n}{(\frac12+\frac13t)_n}
\right|_{n=\frac13t-\frac12}&=&2^{2-t}\cos(\dfrac16\pi t)\label{pocd}
\end{eqnarray}
where $(a)_n:=\Gamma(a+n)/\Gamma(a)$ is the Pochhammer symbol.
They are easily proved, using the reflection and duplication properties
of the gamma function.
In~(\ref{cona},\ref{conb},\ref{conc},\ref{cond})
one sums the Pochhammers of~(\ref{poca}--\ref{pocd})
with weights that are singular for $t/k=n+\frac12$,
with $k=2,1,2,3$, respectively. Using the identity
\begin{equation}
\sum_{n\ge0}\left(\frac{1}{n+\frac12-t/k}\,-\,\frac{1}{n+\frac12+t/k}\right)
=\pi\tan(\pi t/k)\label{tan}
\end{equation}
one then sees that the singularities of the $W$ series are
indeed trigonometric.

Now the crux of these observations is that~(\ref{poca}--\ref{pocd})
are indeed peculiar; there are no further possibilities for
trigonometric reduction of singularities of series of type~(\ref{ww}).
In particular, one readily proves by exhaustion that no
pair of Pochhammers can produce residues at $n=\frac15t-\frac12$
proportional to $2^{-t}\cos(\dfrac1{10}\pi t)$, as required
by~(\ref{what}).
With two pairs of Pochhammers, one has
\begin{equation}
\left.\frac{(\frac12)_n(\frac12-\frac1{10}t)_n}
{\left\{(\frac12+\frac15t)_n\right\}^2}
\right|_{n=\frac15t-\frac12}=2^{3-t}\cos(\dfrac1{10}\pi t)\label{poc4}
\end{equation}
which proves that the Euler-sum generating ${}_4F_3$ series
\begin{equation}
\sum_{n\ge0}\left(\frac{\frac18t}{n+\frac12-\frac15t}
-\frac{\frac18t}{n+\frac12+\frac15t}\right)
\frac{(\frac12)_n(\frac12-\frac1{10}t)_n}
{\left\{(\frac12+\frac15t)_n\right\}^2}
=\dfrac1{40}(\pi t)^2+O(t^3)\label{f43}
\end{equation}
has singularities at positive $2t=5$~mod~10 that are identical to those
of the generator~(\ref{what}). However, the rational values of~(\ref{f43})
at positive $t=5$~mod~10 have denominators with squares of primes
while only single powers occur in the summands of~(\ref{ratp}).
Nor did we succeed in exploiting the more promising identity
\begin{equation}
\left.\frac{\left\{(\frac12-\frac1{10}t)_n\right\}^3}
{\left\{(\frac12)_n\right\}^2(\frac12+\frac15t)_n}
\right|_{n=\frac15t-\frac12}=2^{3-t}\cos(\dfrac1{10}\pi t)\label{poc6}
\end{equation}
which is the simplest way of generating residues with
three pairs of Pochhammer symbols, as might occur in a ${}_5F_4$
representation.

Thus, with reluctance, we leave the reader with a pretty puzzle:
to discover a hypergeometric representation of the generator~(\ref{what})
that reproduces our data, namely
\begin{enumerate}
\item the complex recurrence relation~(\ref{hrec}),
\item the trigonometrically generated singularities of~(\ref{U}),
\item the appearance of Euler sums in expansion~(\ref{expu}),
\item the rational asymptotic expansion~(\ref{uas}),
      generated by the contour integral~(\ref{cont}),
\item the denominator structure of the rational
      values~(\ref{ratp}--\ref{rats}) and~(\ref{rathp},\ref{rathm}).
\end{enumerate}
There is clearly a superabundance of data. Yet until one finds
{\em a priori\/} derivations to replace the existing {\em a posteriori\/}
proofs of~(\ref{cis},\ref{dis}) for the $p=3$
cases~(\ref{leftc},\ref{leftd}), a hypergeometric representation of the
$p=5$ case~(\ref{what}) may remain obscure. When it is found,
the generator of the non-eulerian sums in~(\ref{hn}) will follow,
immediately, from~(\ref{hrec}).

\section{Computation of digits}

Notwithstanding the unsolved puzzle of section~3.5, it is now proven,
by the polylogarithmic analysis of sections~2.1--2.5 and
the hypergeometric analysis of sections~3.1--3.3, that the $d$th
hexadecimal digits of the 18 constants
\begin{eqnarray*}
&&\pi,\phantom{^2}\quad\log2,\\
&&\pi^2,\quad\log^22,\quad\pi\log2,\quad G,\\
&&\pi^3,\quad\log^32,\quad\pi\log^22,\quad\pi^2\log2,\quad\zeta(3),\\
&&\pi^4,\quad\log^42,\quad\pi^2\log^22,\\
&&\phantom{\pi^5,}\quad\log^52,\quad\pi^2\log^32,\quad\pi^4\log2,
\quad\zeta(5),
\end{eqnarray*}
are computable in logarithmic space and almost linear time.
Previously only the first 4 constants were known to be in
this class~\cite{BBP}.

The 64 hexadecimal digits of $\zeta(3)$ that begin at the
10,000,000th place were computed from~(\ref{z3}) to be
$${\tt
CDA018F4E167F435B2AB045FB045A42F86BED12EF82BE2E1C6ECD305E92C5E4B
}\ldots$$
and the corresponding string for $\zeta(5)$ was obtained
from~(\ref{z5}) as $${\tt
F7A15E1277F7B2C04106F04B05C48AC71ACECAB14D555FDA6E5E1EC299535511
}\ldots$$
For $\zeta(5)$ we used David Bailey's
{\sc transmp}~\cite{TRANSMP} to translate 100 lines of
{\sc fortran} to 312--bit precision, sufficient for multiplication modulo
integers up to $62651\times(8\times10^7)^5<2^{148}$,
entailed by the terms in~(\ref{z5}). The 256--bit result
was obtained in 19 hours on a 333~MHz DecAlpha 600 machine,
using merely $0.3$~MB of memory.

It would now be a routine matter to compute the billionth digit of Catalan's
constant, since our result~(\ref{G})
has a simplicity comparable to that of the formul{\ae} for $\pi^2$ and
$\log^22$, computed at this depth in~\cite{BBP}.
More interesting, perhaps, would be a comparison of our results
for $\zeta(3)$ and $\zeta(5)$ with methods based on Wilf--Zeilberger
acceleration~\cite{AZ} of Ap\'ery--like~\cite{BB} results, which
appear~\cite{SP} not yet to have attained the depths probed
by the strings above.

{\em Note added:} In~\cite{D3} we have shown that massive Feynman
diagrams with 3 loops involve SC$^*(2)$ constants,
from the present work, and novel SC$^*(3)$ constants, whose base
of super-fast computation is $b=3$.

{\bf Acknowledgements:} My interest in finding new members of
SC$^*$ came from discussions with Simon Plouffe, during a visit to the
Center for Experimental and Constructive Mathematics at Simon Fraser
University, generously hosted by Jon Borwein in December 1996.
David Bailey of NAS, at NASA--Ames, supplied a finely tuned
version of {\sc pslq} for 16000--bit precision work that indicated
the completeness of integer relations found with Tony Hearn's {\sc reduce},
either analytically or at 128--bit precision. Dirk Kreimer helped me endure
the impasse of section~3.5.

\raggedright

\end{document}